\documentclass[12pt]{amsart}

\usepackage{amsmath,amssymb,amsthm}
\usepackage[all]{xy}
\usepackage{amscd}
\usepackage{amsfonts}
\usepackage{color}
\usepackage[dvipdfm,colorlinks=true]{hyperref}

\numberwithin{equation}{section}

\SelectTips{eu}{12}

\theoremstyle{definition}

\theoremstyle{remark}

\newcommand{\A}{\underline{A}}
\newcommand{\X}{\underline{X}}
\newcommand{\link}{\mathrm{link}}
\newcommand{\MF}{\mathrm{MF}}

\makeatletter
  \def\wideubar{\underaccent{{\cc@style\underline{\mskip10mu}}}}
  \makeatother

\title{On the moment-angle manifold constructed by\\ Fan, Chen, Ma, and Wang}

\author{Kouyemon Iriye}
\address{Department of Mathematics and Information Sciences, Osaka Prefecture University, Sakai, 599-8531, Japan}
\email{kiriye@mi.s.osakafu-u.ac.jp}
\thanks{K.I. is supported by JSPS KAKENHI (No. 26400094).}


\begin{document}

\maketitle

\begin{abstract}
Fan, Chen, Ma, and Wang \cite{FCMW} constructed a moment-angle manifold, $Z_K$, whose cohomology ring 
is isomorphic to that of the connected sum of sphere products consisting of one product of three spheres. In this paper, 
we show that these are in fact diffeomorphic.
\end{abstract}

\baselineskip 16pt

\section{Introduction}

The topology of moment-angle manifolds has been studied by many authors \cite{Mc, DJ, BM, GL}, and it is now 
known that this can be rather complicated. A connected sum of sphere products 
gives a typical example of such manifolds. 

\bigskip

{\bf Theorem 1.1} (McGavran \cite{Mc} and Bosio-Meersseman \cite{BM}){\bf .} 
{\it Let $K$ be the triangulation of a sphere that is dual to the 
simple polytope obtained from the $k$-simplex by cutting off $\ell>0$ vertices, i.e.,
the boundary of a stacked polytope. Then, 
the moment-angle manifold associated to $K$ is diffeomorphic to a connected sum of sphere 
products $Z_K\cong \#_{j=1}^\ell \left(S^{j+2}\times S^{2k+\ell-j-1}\right)^{\#j\binom{\ell+1}{j+1}}.$
}

\bigskip

Here, $X^{\#j}$ denotes the connected sum of $j$-copies of a manifold $X$ without boundary. 
Moreover, for $k=2,3$ Bosio and Meersseman characterized precisely the spherical triangulation that
gives rise to a connected sum of sphere products as a moment-angle manifold. See Proposition 11.6 of 
\cite{BM}. In addition, see the paper \cite{GL}. 
In these observations, only a product of two spheres appears. 

Fan, Chen, Ma, and Wang \cite{FCMW} found that the cohomology ring of the moment-angle manifold 
$Z_{\partial P^8_{28}}$ is isomorphic to that of the connected sum of sphere products
\begin{equation}
M=(S^3\times S^3\times S^6)\#(S^5\times S^7)^{\#8}\# (S^6\times S^6)^{\#8},
\end{equation}
where $\partial P^8_{28}$ is the boundary of $P^8_{28}$ (a $4$-polytope with 18 facets), described in \cite{GS}.  

In this paper, we show that the moment-angle manifold $Z_{\partial P^8_{28}}$ is in fact diffeomorphic 
to the connected sum of sphere products given above, and we affirm the conjecture of \cite{FCMW}. 

\bigskip

{\bf Theorem 1.2.} {\it The moment-angle manifold $Z_{\partial P^8_{28}}$ is diffeomorphic to 
$M$ as defined by $\mathrm{(1.1)}$.}

\bigskip

The author thanks Daisuke Kishimoto for useful conversations, and Editage (www.editage.jp) for English language editing.

\section{Moment-angle manifold}

To prove that $Z_{\partial P^8_{28}}$ is diffeomorphic to $M$, we follow the standard method, i.e., 
the use of $h$-cobordism theory. 
The following theorem is a simple modification of Theorem A1 of Gitler and L\'{o}pez \cite{GL}.  

\bigskip

{\bf Theorem 2.1.} {\it Let $Q$ be a compact smooth manifold of dimension $d+1\geq6$ with 
boundary $\partial Q$, satisfying the following:
\begin{enumerate}
\item[(1)] $Q$ is simply connected with a simply connected boundary.
\item[(2)] There is a finite collection $\{X_j\}$ of disjointly embedded closed smooth manifolds inside $Q$ 
with trivial normal bundles.
\item[(3)] The embedding $\coprod_j X_j\to Q$ induces isomorphisms of integral homology groups of positive 
dimensions, and $H_i(Q)=0$ for $i\geq d-1$. 
\end{enumerate} 
Then, $Q$ is diffeomorphic to the boundary connected sum of 
$
\coprod_j X_j\times D^{d+1-\dim X_j}, 
$
and therefore $\partial Q$ is diffeormophic to 
$
\#_j(X_j\times  S^{d-\dim X_j}). 
$
}

\bigskip

To construct a manifold $Q$ with boundary $Z_{\partial P^8_{28}}$, 
we employ the polyhedral product. 

Let $K$ be a simplicial complex on the vertex set $[m]=\{1,2,\cdots,m\}$. 
A moment-angle complex $Z_K$ associated to $K$ is defined by 
\[
Z_K=\bigcup_{\sigma\in K}(D^2)^\sigma\times(S^1)^{[m]\setminus\sigma}\subset(D^2)^m.
\]
 This construction has been generalized to the polyhedral product. 
See \cite{BBCG, BP}. 
Let $(\X,\A)=\{(X_i,A_i)\}_{i\in[m]}$ be a set of pairs of spaces $(X_i,A_i)$. For a subset $I\subset[m]$, we define 
\[
(\X,\A)^I=\{(x_1,\cdots,x_m)\in \prod_{i=1}^mX_i\ |\ \text{$x_j\in A_j$ if $j\not\in I$}\}.
\] 
Then, the polyhedral product $Z_K(\X,\A)$ is defined 
as the union of $(\X,\A)^\sigma$ for all faces $\sigma$ of $K$:
\[
Z_K(\X,\A)=\bigcup_{\sigma\in K}(\X,\A)^\sigma. 
\]

If $K$ is an $n$-dimensional triangulated sphere with $m$ vertices, 
then $Z_K$ has the structure of an $n+m+1$-dimensional manifold. See Theorem 4.1.4 of \cite{BP}.  
Moreover, if $K$ is the boundary of a simplicial polytope, then $Z_K$ is a smooth manifold. See Theorem 6.2.4 
of \cite{BP}. Therefore, $Z_K$ is called a moment-angle manifold if $K$ is a triangulated sphere. 

Now, we recall the simplicial complex $\partial P^8_{28}$ in the manner described by  
Fan, Chen, Ma, and Wang \cite{FCMW}.  

\bigskip
{\rm Definition 2.2.} For a simplicial complex $K$ on $[m]$, a subset $I \subset[m]$ is called a 
{\it missing face} or {\it minimal non-face} of $K$ if $I$ is not a simplex of $K$, but all of its proper subsets 
are faces of $K$. The set of missing faces of $K$ is denoted by $\MF(K)$. 

\bigskip
For a sequence of integers $i_1<i_2<\cdots<i_k$, $i_1i_2\cdots i_k$ denotes the set $\{i_1, i_2, \cdots, i_k\}$. 
$\partial P^8_{28}$ is characterized by its missing faces
\[
\MF(\partial P^8_{28})=\{56,\ 78,\ 123,\ 134,\ 235,\ 346,\ 147,\ 467,\ 128,\ 258\},
\]
and $\partial P^8_{28}$ has 18 facets (maximal faces):  
\begin{align*}
1245,\ 1246,\ 1257,\ 1267,\ 1357,\ 1367,\ 2347,\ 2367,\ 2457,&\\ 
3457,\ 1458,\ 1468,\ 1358,\ 1368,\ 2348,\ 2368,\ 2468,\ 3458.&
\end{align*}

Let $D^3_+=\{(x,z)\in \mathbf{R}\times\mathbf{C}\ |\ x^2+|z|^2\leq1,\ x\geq0\}$ and 
$S^2_+=\{(x,z)\in D^3_+\ |\ x^2+|z|^2=1\}$. 
Set $Q=Z_{\partial P^8_{28}}((D^3_+,S^2_+),(D^2,S^1),\cdots,(D^2,S^1))$. Then, $Q$ is a 13-dimensional smooth 
manifold with boundary $Z_{\partial P^8_{28}}$. Because $(D^3_+,S^2_+)$ is contractible as a pair of spaces, 
$Q$ is homotopy equivalent to $Z_{\partial P^8_{28}}((*,*),(D^2,S^1),\cdots,(D^2,S^1))$, 
which is identified with $Z_{\partial P^8_{28}-1}$, where $\partial P^8_{28}-1=\{\sigma\in\{2,\cdots,8\}\ |\ \sigma\in \partial P^8_{28}\}$. 
 To apply Theorem 2.1, we need to know the homotopy type of $Z_{\partial P^8_{28}-1}$. 
For a subset $I\subset[8]$, we denote the full subcomplex of $\partial P^8_{28}$ 
on $I$ by $(\partial P^8_{28})_I$. 

\bigskip
{\bf Theorem 2.3.} {\it The following homotopy equivalence holds:
\begin{multline*}
Q\simeq Z_{\partial P^8_{28}-1}\simeq (S^3_{56}\times S^3_{78})\vee S^5_{235}\vee S^5_{346}\vee S^5_{467}\vee  S^5_{258}\\
\vee S^6_{4678}\vee S^6_{4567}\vee S^6_{3467}\vee S^6_{3456}\vee S^6_{2578}
\vee S^6_{2568}\vee S^6_{2358}\vee S^6_{2356}
\vee S^7_{45678}\vee S^7_{25678}\vee S^7_{23568} \vee S^7_{34567},
\end{multline*}
where $S^k_I$ denotes the sphere $S^k$, which appears for the first time 
in $Z_{(\partial P^8_{28})_I}$ as a wedge summand. 
}
\section{Fat wedge filtration}

To prove Theorem 2.3, we make use of the fat wedge filtration introduced by Kishimoto and the author of \cite{IK14}. 

Let $K$ be a simplicial complex on $[m]$. 
For each subset $I\subset[m]$, $K_I$ denotes the full subcomplex on $I$. That is, $K_I=\{\sigma\in K\ | 
\ \sigma\subset I\}$. We regard $Z_{K_I}$ as a subspace of $Z_K$, by identifying it with  
$Z_{K_I}\times \{-1\}^{I^c}$. 
On the other hand, the projection $(D^2)^m\to (D^2)^I$ induces the projection $Z_K\to Z_{K_I}$. 
In particular, $Z_{K_I}$ is a retract of $Z_K$.  

Now, we recall the fat wedge filtration of $Z_K$. For $0\leq i\leq m$, the fat wedge filtration of $Z_K$ is given by
\[
Z_K^i=\{(x_1,\cdots,x_m)\in Z_K\ |\ \text{at least $m-i$ of $x_j$ are $-1$}\},
\]
which induces the following filtration of $Z_K$:
\[
Z_K^0=\{*\}\subset Z_K^1\subset\cdots\subset Z_K^i\subset\cdots\subset Z_K^m=Z_K,
\]
where $*=(-1,\cdots,-1)$. Then, it is easy to see that 
\[
Z_K^i=\bigcup_{I\subset[m],\ |I|\leq i}Z_{K_I}. 
\]

The key property of the fat wedge filtration of the moment-angle complex is the following.

\bigskip

{\bf Theorem 3.1}(Theorem 5.1 of \cite{IK14}){\bf .} {\it 
For $i=1,\dots,m$, $Z_K^i$ is obtained from $Z_K^{i-1}$ 
by attaching a cone to the composition of maps 
$\varphi_{K_I}:\Sigma^i|{K_I}|\to Z_{K_I}^{i-1}\xrightarrow{\mathrm{incl}}Z_K^{i-1}$, 
for each $I\subset[m]$ with $|I|=i$.
}

\bigskip

There exist some classes of simplicial complexes whose attaching maps for $Z_K^i$ are trivial. 
One such complex is the fillable complex. 

\bigskip

{\rm Definition 3.2.} A simplicial complex $K$ is {\it fillable} if there are missing faces $L_1,\cdots, L_r$ of $K$ such that $|K\cup \{L_1,\cdots, L_r\}|$ is contractible.  

\bigskip

{\bf Theorem 3.3}(Theorem 7.2 of \cite{IK14}){\bf .} {\it 
If $K$ is fillable, then the attaching map $\varphi_K$ is null homotopic.
}

\section{Proof of Theorem 2.3}

From this point on, $K$ denotes $\partial P^8_{28}$. Here, we remark that $(K-1)_I=K_I$ 
for a subset $I\subset[8]-1$. 

Because $K-1$ does not have ghost vertices, it is easy to see that $Z_{K-1}^1\simeq*$. 

Step I): $Z_{K-1}^2$. Because by Theorem 3.1 we have a cofiber sequence 
\[
\bigvee_{I\subset[8]-1,\ |I|=2}\Sigma^2|K_I|\to Z_{K-1}^1\simeq* \to Z_{K-1}^2
\]
and $|K_I|$ is not contractible for $I\subset[8]-1$ with $|I|=2$ if and only if $I=56$ or $I=78$, 
we have a homotopy equivalence  $Z_{K-1}^2\simeq S^3_{56}\vee S^3_{78}$. 

Step II):  $Z_{K-1}^3$. $|K_I|$ is not contractible for $I\subset[8]-1$ with $|I|=3$ if and only if $I$ is one of the four missing faces with three vertices in $\MF(K)$, it does not contain the vertex $1$, and in each case $K_I\cong\partial\Delta^2$. Clearly, $\partial \Delta^2$ is fillable, and the attaching maps $\varphi_{K_I}$ are all trivial. Thus, we have that
$
Z_{K-1}^3\simeq S^3_{56}\vee S^3_{78}\vee 
S^5_{235}\vee S^5_{346}\vee S^5_{467}\vee S^5_{258}. 
$

Step III):  $Z_{K-1}^4$. In \cite{FCMW}, the authors classified non-contractible full subcomplexes of $K$ with four vertices. 
There are three possible types. If $I$ is of type A, then $I=5678$ and $Z_{K_I}=S^3\times S^3$. If $I$ is of 
type B or type C, then it is a fillable complex, and the attaching map is trivial. 
Thus, we have the homotopy equivalence 
\begin{multline*}
Z_{K-1}^4\simeq (S^3_{56}\times S^3_{78})\vee S^5_{235}\vee S^5_{346}\vee S^5_{467}\vee  S^5_{258}\\
\vee S^6_{4678}\vee S^6_{4567}\vee S^6_{3467}\vee 
S^6_{3456}\vee S^6_{2578}\vee S^6_{2568}\vee S^6_{2358}\vee S^6_{2356}.
\end{multline*}

Step IV):  $Z_{K-1}^5$. $|K_I|$ is not contractible for $I\subset[8]-1$ with $|I|=5$ if and only if $I$ is 
the compliment of one of the four missing faces with three vertices in $\MF(K)$ that contains the vertex $1$. 
In all cases, it is easy to see that $K_I$ is fillable, and the attaching map is trivial. 

\begin{multline*}
Z_{K-1}^5\simeq (S^3_{56}\times S^3_{78})\vee S^5_{235}\vee S^5_{346}\vee S^5_{467}\vee  S^5_{258}
\vee S^6_{4678}\vee S^6_{4567}\vee S^6_{3467}\vee S^6_{3456}\\
\vee S^6_{2578}\vee S^6_{2568}\vee S^6_{2358}\vee S^6_{2356}
\vee S^7_{45678}\vee S^7_{25678}\vee S^7_{23568}\vee S^7_{34567}.
\end{multline*}

Step V). $Z_{K-1}^6$ and $Z_{K-1}$. 
For $|I|=6$, $|K_I|$ is not contractible if and only if $I$ is the compliment of one of the two missing 
faces with two vertices in $\MF(K)$. That is, $I=123456$ or $I=123478$. In both cases 
$I$ contains the vertex $1$, and therefore $|(K-1)_I|$ is contractible for all $I\subset[8]-1$ with $|I|=6$.   

Because $|K-1|$ is contractible, $Z_{K-1}^7\simeq Z_{K-1}^6\simeq Z_{K-1}^5$, and we have completed the proof of 
Theorem 2.3.

\section{Proof of Theorem 1.2}

We call $S^3_{56}\times S^3_{78}$ or $S^k_I$ a wedge summand when they appear 
as a wedge summand of $Q$. To construct an embedding for a wedge summand to $Q$, we make use of 
the smooth embedding $Z_{\link_K(\sigma)}\to Z_K$, where $\link_K(\sigma)=\{\tau\in \sigma^c\ | 
\tau\cup\sigma\in K\}$ is the link of $\sigma$ in $K$. To embed $S^3_{56}\times S^3_{78}$, we use 
$\link_K(13)=\{\emptyset, 5,6\}*\{\emptyset, 7,8\}$. Then, $Z_{\link_K(13)}=S^3\times S^3$ and 
the embedding $Z_{\link_K(13)}=S^3\times S^3\to Z_K\subset Q$ represents the wedge summand 
$S^3_{56}\times S^3_{78}$. 

To embed the other wedge summands, we use $\link_K(2)$ and $\link_K(4)$. 

The facets of $\link_K(2)$ are $145,\ 146,\ 157,\ 167,\ 347,\ 367,\ 457, \ 348,\ 368,\ 468$,  
and therefore it is the boundary of a stacked polytope. By Theorem 1.1, we see that
$Z_{\link_K(2)}\cong(S^3\times S^7)^{\#6}\#(S^4\times S^6)^{\#8}\#(S^5\times S^5)^{\#3}$. 
It is easy to see where a sphere factor of $Z_{\link_K(2)}$ occurs, and we see that 
\begin{eqnarray*}
Z_{\link_K(2)}&\cong&
 (S^3_{13}\times S^7_{45678})\#(S^3_{18}\times S^7_{34567})\#(S^3_{35}\times S^7_{14678})\\
&&\#(S^3_{56}\times S^7_{13478})\#(S^3_{58}\times S^7_{13467})\#(S^3_{78}\times S^7_{13456})\\
&&\#(S^4_{135}\times S^6_{4678})\#(S^4_{138}\times S^6_{4567})\#(S^4_{158}\times S^6_{3467})\#(S^4_{178}\times S^6_{3456})\\
&&\#(S^4_{356}\times S^6_{1478})\#(S^4_{358}\times S^6_{1467})\#(S^4_{568}\times S^6_{1347})\#(S^4_{578}\times S^6_{1346})\\
&&\#(S^5_{147}\times S^5_{3568})\#(S^5_{346}\times S^5_{1578})\#(S^5_{467}\times S^5_{1358}). 
\end{eqnarray*}
Here, we remark that in the formula above, factors of $S^k_I$ arise from $Z_{(\link_K(2))_I}$. 
The composite of maps $Z_{\link_K(2)}\to Z_K\subset Q\xrightarrow{\simeq} Z_{K-1}$ maps a summand $S^k_I$ 
in $Z_{\link_K(2)}$ to $S^k_I$ in $Z_{K-1}$ in a homotopy equivalent manner if $1\not\in I$, and $*$ otherwise.  
To see this, it is sufficient to check that $(\link_K(2))_I\to K_I$ is homotopy equivalent if $S^k_I$ appears 
as a factor of the sphere products in $Z_{\link_K(2)}$ and $1\not\in I$. 
It is easy to verify this fact in all cases, and so we omit the details.

Similarly, the facets of $\link_K(4)$ are 
$
126,\ 125,\ 168, \ 158,\ 238,\ 358,\ 268,\ 237,\ 357,\ 257.
$
These are just those of $\link_K(2)$, except with $4$ replaced by $2$ and $7\leftrightarrow8$ exchanged with
$5\leftrightarrow6$. Therefore, we see that 
\begin{eqnarray*}
Z_{\link_K(4)}&\cong&
 (S^3_{13}\times S^7_{25678})\#(S^3_{17}\times S^7_{23568})\#(S^3_{36}\times S^7_{12578})\\
&&\#(S^3_{56}\times S^7_{12378})\#(S^3_{67}\times S^7_{12358})\#(S^3_{78}\times S^7_{12356})\\
&&\#(S^4_{136}\times S^6_{2578})\#(S^4_{137}\times S^6_{2568})\#(S^4_{167}\times S^6_{2358})\#(S^4_{178}\times S^6_{2356})\\
&&\#(S^4_{356}\times S^6_{1278})\#(S^4_{367}\times S^6_{1258})\#(S^4_{567}\times S^6_{1238})\#(S^4_{678}\times S^6_{1235})\\
&&\#(S^5_{128}\times S^5_{3567})\#(S^5_{235}\times S^5_{1678})\#(S^5_{258}\times S^5_{1367}). 
\end{eqnarray*}

Thus, we have obtained a necessary embedding for each wedge summand. Because all wedge summands 
are embedded in the boundary of $Q$, the embeddings can be made mutually disjoint and
lying in a collar neighborhood of the boundary $Z_K$.   

Finally, we see that the normal bundles of the embedded manifolds are trivial. The normal bundle of $S^5$ 
(resp. $S^7$) embedded in $Z_K$ has been classified by the homotopy group $[S^5,BO(7)]\cong[S^4,O(7)]\cong0$ 
(resp.  $[S^7,BO(5)]\cong[S^6,O(5)]\cong0$) in \cite{MT} and \cite{MM}. Thus, these are also trivial in $Q$.  
The normal bundle of $S^6$ embedded in $Q$ is classified by the homotopy group 
$[S^6,BO(7)]\cong[S^5,O(7)]\cong0$. The normal bundle of $S^3\times S^3$ embedded in $Q$ is classified 
by the homotopy set $[S^3\times S^3, BO(7)]=*$. To see this, we consider the exact sequence associated with the 
cofiber sequence $S^3\vee S^3\to S^3\times S^3\to S^6$:
\[
0\cong[S^6, BO(7)]\to [S^3\times S^3, BO(7)]\to [S^3\vee S^3, BO(7)].
\]
Because $[S^3\vee S^3, BO(7)]\cong[S^3,BO(7)]\times [S^3,BO(7)]\cong [S^2,O(7)]\times [S^2,O(7)]\cong0$, 
we see that $ [S^3\times S^3, BO(7)]=*$. Thus, all of the normal bundles are trivial, and we have proved 
Theorem 1.2 by applying Theorem 2.1.

\end{document}